\documentclass[11pt]{article}
\usepackage{makeidx}
\usepackage{amsmath}
\usepackage{amsfonts}
\usepackage{amssymb}
\usepackage{graphicx}
\usepackage{multicol}

\DeclareMathOperator{\dom}{dom}

\DeclareMathOperator{\gph}{gph}
\DeclareMathOperator{\clm}{clm}
\DeclareMathOperator{\Lipusc}{Lipusc}

\DeclareMathOperator{\extr}{extr}
\DeclareMathOperator{\spann}{span}

\DeclareMathOperator{\conv}{conv}

\DeclareMathOperator{\bd}{bd}
\DeclareMathOperator{\ri}{ri}
\DeclareMathOperator{\dist}{d\!}

\newtheorem{theo}{Theorem}
\newtheorem{lem}{Lemma}
\newtheorem{prop}{Proposition}
\newtheorem{cor}{Corollary}
\newtheorem{rem}{Remark}
\newtheorem{exa}{Example}

\newenvironment{dem}[1][Proof]{\noindent \textbf{#1.} }{\ \rule{0.5em}{0.5em}}
\begin{document}

\title{Lipschitz upper semicontinuity of linear inequality systems under
full perturbations \thanks{%
This research has been partially supported by Grant PID2022-136399NB-C22
from MICINN, Spain, and ERDF, "A way to make Europe", European Union.}}
\author{J. Camacho\thanks{%
Center of Operations Research, Miguel Hern\'{a}ndez University of Elche,
03202 Elche (Alicante), Spain (j.camacho@umh.es, canovas@umh.es,
parra@umh.es).} \and M.J. C\'{a}novas\footnotemark[2] \and H. Gfrerer\thanks{%
Johann Radon Institute for Computational and Applied Mathematics (RICAM),
A-4040, Linz, Austria (helmut.gfrerer@ricam.oeaw.ac.at).} \and J. Parra\footnotemark[2]%
\and \small{Dedicated to Michel Théra on the occasion of his 80th birthday}}
\date{}
\maketitle

\begin{abstract}
The present paper is focused on the computation of the Lipschitz upper
semicontinuity modulus of the feasible set mapping in the context of fully
perturbed linear inequality systems; i.e., where all coefficients are
allowed to be perturbed. The direct antecedent comes from the framework of
right-hand side (RHS, for short) perturbations. The difference between both
parametric contexts, full \emph{vs }RHS perturbations, is emphasized. In
particular, the polyhedral structure of the graph of the feasible set
mapping in the latter framework enables us to apply classical results as
those of Hoffman \cite{Hof52} and Robinson \cite{Rob81}. In contrast, the
graph of the feasible set mapping under full perturbations is no longer
polyhedral (not even convex). This fact requires \emph{ad hoc }techniques
to analyze the Lipschitz upper semicontinuity property and its corresponding
modulus.\bigskip


\bigskip


\end{abstract}

\section{Introduction}

This paper is focused on the \emph{upper} Lipschitzian behavior of the
feasible set associated with a parameterized linear inequality system
written in the form%
\begin{equation}
\sigma \left( A,b\right) :=\left\{ a_{t}^{\prime }x\leq b_{t},~t\in
T:=\{1,...,m\right\} \},  \label{sigma(a,b)}
\end{equation}%
where $x\in \mathbb{R}^{n}$ is the vector of variables, $\left( A,b\right)
\in \mathbb{R}^{m\times n}\times \mathbb{R}^{m}$ is the parameter to be
perturbed, $a_{t}^{\prime }$ represents the $t$-th row of matrix $A,$ and $%
b_{t}$ is the $t$-th coordinate of vector $b.$ Vectors in $\mathbb{R}^{n}$
are understood as column-vectors and the prime represents transposition.
Associated with the parameterized family of systems (\ref{sigma(a,b)}) we
consider the \emph{feasible set mapping }$\mathcal{F}:\mathbb{R}^{m\times
n}\times \mathbb{R}^{m}\rightrightarrows \mathbb{R}^{n}$\emph{\ }given by
\begin{equation}
\mathcal{F}\left( A,b\right) :=\left\{ x\in \mathbb{R}^{n}\mid Ax\leq
b\right\} .  \label{eq_F}
\end{equation}%
Given a fixed $\bar{A}\in \mathbb{R}^{m\times n},$ we also consider the
multifunction $\mathcal{F}_{\bar{A}}:\mathbb{R}^{m}\rightrightarrows \mathbb{%
R}^{n}$ defined as%
\begin{equation*}
\mathcal{F}_{\bar{A}}\left( b\right) :=\mathcal{F}\left( \bar{A},b\right) ,%
\text{ for all }b\in \mathbb{R}^{m}.
\end{equation*}%
In this way, $\mathcal{F}_{\bar{A}}$ is set in the framework of \emph{%
right-hand side} (RHS, in brief) perturbations, while $\mathcal{F}$
formalizes the variation of the feasible set under \emph{full }perturbations.

The main goal of the paper is to analyze the \emph{Lipschitz upper
semicontinuity }property of multifunction $\mathcal{F}$ at a nominal (given)
element $\left( \bar{A},\overline{b}\right) \in \mathrm{dom}\mathcal{F}$
(where `$\mathrm{dom}$' stands for domain; i.e., $\left( \bar{A},\overline{b}%
\right) \in \mathrm{dom}\mathcal{F\Leftrightarrow F}\left( \bar{A},\overline{%
b}\right) \neq \emptyset $) as well as to compute the corresponding modulus
of $\mathcal{F}$ at $\left( \bar{A},\overline{b}\right) ,$ denoted by $%
\Lipusc\mathcal{F}\left( \bar{A},\overline{b}\right) .$ We are appealing to
the terminology of Klatte and Kummer \cite{KlKu02} or Uderzo \cite{Uderzo21}%
, although this property was introduced by Robinson \cite{Rob81} under the
name of \emph{upper Lipschitz continuity }and can be found in the monograph
of Dontchev and Rockafellar \cite{DoRo} as\emph{\ outer Lipschitz continuity}%
. The immediate antecedent to the present work can be traced out from \cite%
{CCP21}, where a formula for computing the Lipschitz upper semicontinuity
modulus of $\mathcal{F}_{\bar{A}}$ at a nominal $\overline{b}\in \mathrm{dom}%
\mathcal{F}_{_{\bar{A}}},$ $\Lipusc\mathcal{F}_{\bar{A}}\left( \overline{b}%
\right) ,$ is provided. Roughly speaking, these moduli constitute the
`tightest' measures of the expansion of feasible sets with respect to data
perturbations, around a nominal $\left( \bar{A},\overline{b}\right) \in
\mathrm{dom}\mathcal{F}$ or $\overline{b}\in \mathrm{dom}\mathcal{F}_{_{\bar{%
A}}},$ depending on the context; see the following paragraphs for some
preliminary details and Section 2 for the definitions and first results
about generic multifunctions between metric spaces.

At this moment, we point out the fact that there are notable differences
between both contexts of RHS and full perturbations. To start with, $\mathrm{%
gph}\mathcal{F}_{\bar{A}}$ (where `$\mathrm{gph}$' means graph) is a convex
polyhedral cone as the solution set of system `$-b+Ax\leq 0_{m}$' in the
variable $\left( b,x\right) \in \mathbb{R}^{m}\times \mathbb{R}^{n}.$ Hence
appealing to a classical result by Robinson \cite{Rob81} (which applies for
mappings whose graphs are finite union of convex polyhedra), $\mathcal{F}_{%
\bar{A}}$ is Lipschitz upper semicontinuous at any element $\overline{b}\in
\mathrm{dom}\mathcal{F}_{_{\bar{A}}}$, which reads as the existence of a
constant $\kappa \geq 0$ along with a neighborhood $V\subset \mathbb{R}^{m}$
of $\overline{b}$ such that
\begin{equation}
\dist\left( x,\mathcal{F}_{\bar{A}}\left( \overline{b}\right) \right) \leq
\kappa \dist\left( b,\overline{b}\right) \text{ for all }b\in V\text{ and
all }x\in \mathcal{F}_{_{\bar{A}}}\left( b\right) ,  \label{eq_Lipusc_F_A}
\end{equation}%
where `$\dist$' denotes the distance in both spaces of variables and
parameters. By definition, $\Lipusc\mathcal{F}_{\bar{A}}\left( \overline{b}%
\right) $ is the infimum of all constants $\kappa \geq 0$ verifying (\ref%
{eq_Lipusc_F_A}). Hence $\Lipusc\mathcal{F}_{\bar{A}}\left( \overline{b}%
\right) $ is finite for all $\overline{b}\in \mathrm{dom}\mathcal{F}_{_{\bar{%
A}}}.$ However, $\mathrm{gph}\mathcal{F}$ is no longer convex (and not a
cone) which leads to a different scenario. Indeed, we advance that $\Lipusc%
\mathcal{F}\left( \bar{A},\overline{b}\right) $ is finite if and only if one
of the following conditions holds: $\left( i\right) $ $\mathcal{F}\left(
\bar{A},\overline{b}\right) $ is bounded, $\left( ii\right) $ $\mathcal{F}%
\left( \bar{A},\overline{b}\right) =\mathbb{R}^{n},$ or $\left( iii\right) $
$n=1$ (see Proposition \ref{Prop_F(A,b)_unbounded} and Theorem \ref{theo max
clm} below).

Looking at (\ref{eq_Lipusc_F_A}), the Lipschitz upper semicontinuity
property can be seen as a \emph{semilocal} Lipschitzian property in the
sense that it involves the whole image set (all points $x\in \mathcal{F}_{_{%
\bar{A}}}\left( b\right) )$ but only parameters around the nominal $%
\overline{b}$ (recall $b\in V).$ For completeness, and with the aim of
emphasizing the difference between the contexts of RHS and full
perturbations, let us comment what happens with the so-called \emph{Hoffman
constant. }The celebrated Hoffman lemma (see \cite{Hof52}) entails the
existence of $\kappa \geq 0$ such that%
\begin{equation}
\dist\left( x_{1},\mathcal{F}_{\bar{A}}\left( b_{2}\right) \right) \leq
\kappa \dist\left( b_{1},b_{2}\right) \text{ for all }b_{1},b_{2}\in \mathrm{%
dom}\mathcal{F}_{_{\bar{A}}}\text{ and all }x_{1}\in \mathcal{F}_{_{\bar{A}%
}}\left( b_{1}\right) .  \label{eq_Hof_lip}
\end{equation}%
This last property, which is known as the \emph{Lipschitz continuity} of $%
\mathcal{F}_{_{\bar{A}}}$ \emph{on its domain}, is obviously stronger than
the Lipschitz upper semicontinuity. Observe that condition (\ref{eq_Hof_lip})%
\emph{\ }provides a \emph{global} (instead of semilocal) Lipschitzian
property as it involves all parameters in the whole $\mathrm{dom}\mathcal{F}%
_{_{\bar{A}}}$ (instead of parameters around a fixed $\overline{b}).\,$The
(sharp) Hoffman constant of $\mathcal{F}_{_{\bar{A}}},$ denoted by $\mathrm{%
Hof}\mathcal{F}_{_{\bar{A}}},$ is the infimum of all $\kappa \geq 0\,\ $%
verifying (\ref{eq_Hof_lip}). Constant $\mathrm{Hof}\mathcal{F}_{_{\bar{A}}}$
is a relevant measure of the stability of $\mathcal{F}_{\bar{A}}$ and
different authors have provided computable expressions for it; let us cite Az%
\'{e} and Corvellec \cite[Theorem 3.3]{AzCo02}, Burke and Tseng \cite[%
Theorem 8]{BuTseng96}, Klatte and Thiere \cite[Theorem 2.7]{KlTh95}, Li \cite%
{Li93}, Pe\~{n}a \emph{et al. }\cite[Formula (3)]{PVZ20}, and Z\u{a}linescu
\cite{Za03}, among others. However, the counterpart of (\ref{eq_Hof_lip}) in
the context of full perturbations is no longer true and the Hoffman constant
of $\mathcal{F}$ is infinite. In fact, in the case $n\geq 2,$ just taking
any $\left( \bar{A},\overline{b}\right) \in \mathrm{dom}\mathcal{F}$ such
that $\mathcal{F}\left( \bar{A},\overline{b}\right) $ is unbounded and
different from $\mathbb{R}^{n}$, it is clear from the definitions that this
Hoffman constant is greater than or equal to $\Lipusc\mathcal{F}\left( \bar{A%
},\overline{b}\right) =+\infty .$ In the case $n=1,$ for the system
consisting of the single inequality $x\leq \alpha ,$ letting $\alpha
\rightarrow +\infty ,$ Proposition \ref{Prop_F(A,b)_unbounded} below also
ensures that the Hoffman constant of $\mathcal{F}$ is infinite in this case.

Apart from the announced differences between the stability of $\mathcal{F}$
and $\mathcal{F}_{_{\bar{A}}},$ let us advance that, via \emph{ad hoc }%
techniques, this paper brings to light a remarkable similarity.
Specifically, Theorem \ref{theo max clm}, together with Corollary \ref%
{Cor_Lipusc_F}, provides a point-based formula (based on the nominal data)
for $\Lipusc\mathcal{F}\left( \bar{A},\overline{b}\right) $ when $\mathcal{F}%
\left( \bar{A},\overline{b}\right) $ is bounded, which reminds of the one
given in \cite[Theorem 5]{CCP21} for $\Lipusc\mathcal{F}_{\bar{A}}\left(
\overline{b}\right) $ (recalled here in Theorem \ref{theo Lipusc vs clm
finite}). In both cases, the Lipschitz upper semicontinuity moduli are
expressed in terms of the supremum of certain \emph{local }stability
measures (called \emph{calmness }moduli) evaluated in the set of extreme
points of the nominal feasible set.

Regarding the structure of the paper, Section 2 recalls the necessary
concepts and previous results which are appealed to in the current paper.
Section 3 contains the main results of this work, devoted to the computation
of $\Lipusc\mathcal{F}\left( \bar{A},\overline{b}\right) ,$ with $\left(
\bar{A},\overline{b}\right) \in \mathrm{dom}\mathcal{F},$ by distinguishing
two cases focused on the boundedness or not of $\mathcal{F}\left( \bar{A},%
\overline{b}\right) .$ Finally, Section 4 presents some conclusions,
illustrations and comments on further research.

\section{Preliminaries}

\label{sec: pre} Let us start by fixing the notation and introducing the
main definitions used throughout the paper. Given $X\subset \mathbb{R}^{p}$,
$p\in \mathbb{N}$, we denote by $\mathrm{int}X$, $\mathrm{ri}X$, $\mathrm{%
conv}X$, and $\mathrm{span}X$ the interior, the relative interior, the
convex hull, and the linear hull of $X$ respectively. Hereafter, the
convention $\mathrm{conv}\emptyset =\emptyset $ is adopted. Moreover, if $X$
is convex, its set of extreme points is denoted by $\mathrm{extr}X$.

We appeal to several Lipschitz type properties that are recalled now for a
generic multifunction between metric spaces; the reader is addressed to the
monographs \cite{DoRo, KlKu02, mor06a, rw} for extra information about
Lipschitz type properties among other topics of variational analysis. Let $%
\mathcal{M}:Y\rightrightarrows X$ be a set-valued mapping between metric
spaces with both distances denoted by $\dist$ indistinctly. This paper
focuses on the Lipschitz upper semicontinuity property of $\mathcal{M}$ at $%
\overline{y}\in \dom\mathcal{M},$ which is defined as the existence of a
neighborhood $V$ of $\overline{y}$ together with a constant $\kappa \geq 0$
such that
\begin{equation}
\dist\left( x,\mathcal{M}\left( \overline{y}\right) \right) \leq \kappa \dist%
\left( y,\overline{y}\right) \text{ for all }y\in V\text{ and all }x\in
\mathcal{M}\left( y\right) ,  \label{eq Lipusc property}
\end{equation}%
where $\dist\left( x,\Omega \right) :=\inf \left\{ \dist\left( x,\omega
\right) \mid \omega \in \Omega \right\} $, for $x\in X$ and $\Omega \subset X
$. For convenience we assume $\inf \emptyset :=+\infty ;$ in particular, $%
\dist\left( x,\emptyset \right) =+\infty .$ As commented in the previous
section, the Lipschitz upper semicontinuity property is considered to be a
semilocal measure in the sense that only $y$'s are required to be around the
nominal value, while $x$ may vary in the whole image set $\mathcal{M}\left(
y\right) $.

Concerning local measures, the \emph{calmness} property holds for $\mathcal{M%
}$ at $\left( \overline{y},\overline{x}\right) \in \mathrm{gph}\mathcal{M}$
if there exist a constant $\kappa \geq 0$ and a neighborhood of $\left(
\overline{y},\overline{x}\right) $, $V\times U\subset Y\times X$, such that
\begin{equation}
\dist\left( x,\mathcal{M}\left( \overline{y}\right) \right) \leq \kappa \dist%
\left( y,\overline{y}\right) \text{ for all }x\in \mathcal{M}(y)\cap U\text{
and all }y\in V.  \label{eq calmness property}
\end{equation}%
It can be equivalently read in terms of the inverse function, $\mathcal{M}%
^{-1}$, i.e. $y\in \mathcal{M}^{-1}(x)\iff x\in \mathcal{M}(y)$, appealing
to the well known \emph{metric subregularity} property (cf. \cite[Theorem
3H.3 and Exercise 3H.4]{DoRo}) of $\mathcal{M}^{-1}$ at $\left( \overline{x},%
\overline{y}\right) $. Namely, the existence of $\kappa \geq 0$ and a
(possibly smaller) neighborhood $U$ of $\overline{x}$ such that
\begin{equation}
\dist\left( x,\mathcal{M}\left( \overline{y}\right) \right) \leq \kappa \dist%
\left( \overline{y},\mathcal{M}^{-1}(x)\right) \text{ for all }x\in U.
\label{eq met subreg}
\end{equation}

The infimum of constants $\kappa $ appearing in \eqref{eq Lipusc property}
and \eqref{eq calmness property}, or equivalently \eqref{eq met subreg}, for
some associated neighborhoods, are the corresponding moduli of those
properties. Specifically, $\Lipusc\mathcal{M}\left( \overline{y}\right) $ is
the \emph{Lipschitz upper semicontinuity modulus} of $\mathcal{M}$ at a
given $\overline{y}$ and $\clm\mathcal{M}\left( \overline{y},\overline{x}%
\right) $ denotes the \emph{calmness modulus} of $\mathcal{M}$ at $\left(
\overline{y},\overline{x}\right) $. These constants may be rewritten using
upper limits instead of infima, as the following expressions show:
\begin{flalign}
\Lipusc\mathcal{M}(\overline{y})&=\limsup\limits_{y\rightarrow \overline{y}%
}\left( \sup\limits_{x\in \mathcal{M}(y)}\dfrac{\dist\left(x,\mathcal{M}\left(\overline{y}%
\right)\right)}{\dist\left( y,\overline{y}\right) }\right) ,\text{ }\overline{y}\in \dom%
\mathcal{M},\label{eq Lipusc equiv} \\
\clm\mathcal{M}\left( \overline{y},\overline{x}\right)
&=\limsup\limits_{(y,x)\overset{\gph\mathcal{M}}{\longrightarrow }\left( \overline{y},\overline{x}\right)}\dfrac{\dist\left(x,\mathcal{M}\left(\overline{y}%
\right)\right)}{\dist\left( y,\overline{y}\right) },\text{ }\left(
\overline{y},\overline{x}\right) \in \gph\mathcal{M}, \nonumber
\end{flalign}where $\left( y,x\right) \overset{\gph\mathcal{M}}{%
\longrightarrow }\left( \overline{y},\overline{x}\right) $ means $\left(
y,x\right) \rightarrow \left( \overline{y},\overline{x}\right) $ with $%
\left( y,x\right) \in \gph\mathcal{M}.$ The first equality, (\ref{eq Lipusc
equiv}), was established in \cite[Proposition 2]{CCP21}, while the second
one comes directly from the definition. Here, $\limsup_{z\rightarrow
\overline{z}}$ represents the supremum (maximum, indeed) over all possible
sequences $\left\{ z_{r}\right\} _{r\in \mathbb{N}}\rightarrow \overline{z}$
under the convention $\frac{0}{0}:=0$. Given $\overline{y}\in \dom\mathcal{M}
$, it is straightforward to check that
\begin{equation}
\Lipusc\mathcal{M}(\overline{y})\geq \sup_{x\in \mathcal{M}(\overline{y})}%
\clm\mathcal{M}(\overline{y},x)  \label{eq Lipusc vs clm}
\end{equation}%
General multifunctions may present strict inequality in
\eqref{eq Lipusc vs
clm}, as illustrated by the following example.

\begin{exa}
\emph{\cite[Example 2]{CCP21} Let }$\mathcal{M}:\mathbb{R}\rightarrow
\mathbb{R}$\emph{\ be defined by }$\mathcal{M}(y):=0$\emph{\ for }$y\leq 0$%
\emph{\ and }$\mathcal{M}(y):=1$\emph{\ for} $y>0$\emph{. Then }$\clm%
\mathcal{M}(y,x)=0$\emph{\ for all }$(y,x)\in \gph\mathcal{M}$\emph{, while }%
$\Lipusc\mathcal{M}(0)=+\infty $\emph{\ by just considering the sequence }$%
y_{r}:=\frac{1}{r}$.
\end{exa}

At this point, let us recall an important class of mappings that do satisfy
the equality in \eqref{eq Lipusc vs clm}, as stated in \cite[Theorem 4]%
{CCP21}. We recover this result in the following theorem for completeness.

\begin{theo}
\label{theo Lipusc vs clm convex} Let $\mathcal{M}:Y\rightrightarrows X,$
with $Y$ being a normed space and $X$ being a reflexive Banach space. Assume
that $\gph\mathcal{M}$ is a nonempty convex set. For any $\overline{y}\in %
\dom\mathcal{M}$ such that $\mathcal{M}(\overline{y})$ is closed, we have
\begin{equation}
\Lipusc\mathcal{M}(\overline{y})=\sup_{x\in \mathcal{M}(\overline{y})}\clm%
\mathcal{M}(\overline{y},x).  \label{eq_44}
\end{equation}
\end{theo}

Indeed, the feasible set mapping of linear - possibly infinite - inequality
systems with respect to right-hand side perturbations fits into this class.
Moreover, certain subclasses of mappings attain the equality (\ref{eq_44})
as a maximum in some specific points. In fact, whenever the system is
confined to finitely many constraints, one may select an extension of the
set of extreme points of the original feasible set: (see \cite[p. 142]{Li94}
and \cite[Section 2.2]{GCPT18} for an in-depth discussion on this extension)
\begin{equation}
\mathcal{E}\left( A,b\right) :=\extr\left( \mathcal{F}\left( b\right) \cap %
\spann\left\{ a_{t},~t\in T\right\} \right) ,\;\left( A,b\right) \in \dom%
\mathcal{F},  \label{eq E}
\end{equation}%
and the following theorem shows how we may confine to this finite subset of
points.

\begin{theo}
\emph{(cf. \cite[Theorem 5]{CCP21})} \label{theo Lipusc vs clm finite} Let $%
\left( \overline{b},\overline{x}\right) \in \gph\mathcal{F}_{\bar{A}}$, then
\begin{equation*}
\Lipusc\mathcal{F}_{\bar{A}}\left( \overline{b}\right) =\max_{x\in \mathcal{E%
}\left( \bar{A},\overline{b}\right) }\clm\mathcal{F}_{\bar{A}}\left(
\overline{b},x\right) .
\end{equation*}
\end{theo}

The aim of this paper is to look for the counterpart results of the previous
theorem in the context of full perturbations, where the feasible set mapping
presents a nonconvex graph. Theorem \ref{theo sup clm} fulfills this goal,
providing a new context in which equality in \eqref{eq Lipusc vs clm} holds.
In a second step, the supremum is replaced by a maximum over appropriate
sets of extreme points, captured in Theorem \ref{theo max clm}. Indeed, the
ability to compute the Lipschitz upper semicontinuity modulus is translated
into the possibility of determining the calmness modulus at a fixed point of
the graph of the multifunction (see Section 4 for some illustrations).

The second half of this section provides some background on the calmness
moduli for multifunction $\mathcal{F}$. First, let us fix the topology of
the corresponding spaces. The space of variables, $\mathbb{R}^{n}$, is
endowed with an arbitrary norm $\left\Vert \cdot \right\Vert $, whose dual
norm, $\left\Vert \cdot \right\Vert _{\ast }$, is given by $\left\Vert
u\right\Vert _{\ast }=\max_{\left\Vert x\right\Vert \leq 1}\left\vert
u^{\prime }x\right\vert $. Elements $b$ in the right-hand side parameter
space, $\mathbb{R}^{m}$, work with the norm $\left\Vert b\right\Vert
_{\infty }:=\max_{t\in T}\left\vert b_{t}\right\vert $. The full parameter
space, $\mathbb{R}^{m\times n}\times \mathbb{R}^{m}$, is endowed with the
norm $\left\Vert \left( A,b\right) \right\Vert :=\max_{t\in T}\Vert
(a_{t},b_{t})\Vert $, where $\Vert (a_{t},b_{t})\Vert =\max \left\{
\left\Vert a_{t}\right\Vert _{\ast },\left\vert b_{t}\right\vert \right\} $.
Here $a_{t}$ is identified with the functional $x\mapsto a_{t}^{\prime }x,$
hence the use of the dual norm.

In relation to the feasible set mapping $\mathcal{F}$, the set of active
indices at $x\in \mathcal{F}\left( A,b\right) $ is denoted by
\begin{equation}
T_{A,b}\left( x\right) :=\left\{ t\in T\mid a_{t}^{\prime }x=b_{t}\right\} .
\label{eq T}
\end{equation}%
Furthermore, we appeal to the family of sets $\mathcal{D}_{A,b}(x)$
(introduced in \cite{CLPT14}) conformed by all $D\subset T_{A,b}\left(
x\right) $ such that the following system in the variable $d\in \mathbb{R}%
^{n}$ is consistent:
\begin{equation}
\left\{
\begin{tabular}{rl}
$a_{t}^{\prime }d=1,$ & $t\in D,$ \\
$a_{t}^{\prime }d<1,$ & $t\in T_{A,b}\left( x\right) \setminus D$%
\end{tabular}%
\ \right\} .  \label{eq system D}
\end{equation}%
Note that, for any solution $d\in \mathbb{R}^{n}$ of \eqref{eq system D}
with associated index set $D$, the hyperplane $\left\{ z\in \mathbb{R}%
^{n}\mid d^{\prime }z=1\right\} $ contains the set $\left\{ a_{t},~t\in
D\right\} $ and leaves $\{0_{n}\}\cup \left\{ a_{t},~t\in T\left( x\right)
\setminus D\right\} $ on the same open half-space.

To finish this section, we introduce the main preliminary results that will
be used in the sequel. The first one, which deals with the set of parameters
that make a given point feasible, simply adapts \cite[Lemma 1]{CLPT14} to the
notation and framework of this paper.

\begin{lem}
\label{lem dist inverse} Let $\left( \bar{A},\overline{b}\right) \in \dom%
\mathcal{F}$. For any $x\in \mathbb{R}^{n}$ we have
\begin{equation*}
\dist\left( \left( \bar{A},\overline{b}\right) ,\mathcal{F}^{-1}(x)\right) =%
\frac{\left\Vert \left( \bar{A}x-\overline{b}\right) _{+}\right\Vert
_{\infty }}{\left\Vert x\right\Vert +1},
\end{equation*}%
where $\left( \bar{A}x-\overline{b}\right) _{+}$ denotes the vector with
components $\left( \overline{a}_{t}^{\prime }x-\overline{b}_{t}\right) _{+}$%
, the nonnegative part of $\overline{a}_{t}^{\prime }x-\overline{b}_{t}$, for $%
t\in T$.
\end{lem}

\begin{rem}
\label{Rem 1b}\emph{The reader can easily check that}
\begin{equation*}
\left\Vert \left( \bar{A}x-\overline{b}\right) _{+}\right\Vert _{\infty }=%
\dist\left( \overline{b},\mathcal{F}_{\bar{A}}^{-1}(x)\right)
\end{equation*}%
\emph{for all }$\overline{b}\in \dom\mathcal{F}_{\bar{A}}$\emph{\ and all }$%
x\in \mathbb{R}^{n}.$
\end{rem}

Thiese expressions for the respective distances are considered to be
point-based in the sense that they only depend on the nominal parameter $%
\left( \bar{A},\overline{b}\right) $ and a fixed $x\in \mathbb{R}^{n}$. In
the same spirit, the next result provides a computable formula for the
calmness modulus of the feasible set mapping at a given point. In it, $%
\mathrm{end\,}C$ denotes the end set of the convex set $C\subset \mathbb{R}%
^{n}$ (see \cite{Hu05})\emph{\ }defined as
\begin{equation*}
\mathrm{end\,}C:=\left\{ u\in \mathrm{cl\,}C\mid \nexists \mu >1\text{ such
that }\mu u\in \mathrm{cl\,}C\right\} .
\end{equation*}%
In particular, $\mathrm{end\,}\emptyset =\emptyset $ and $0_{n}\notin
\mathrm{end\,}C$ for any convex subset $C\subset \mathbb{R}^{n}.$

\begin{theo}
\emph{\cite[Theorems 4 and 5]{CLPT14}} \label{theo clm feasible} Let $\left(
\left( \bar{A},\overline{b}\right) ,\overline{x}\right) \in \gph\mathcal{F}$%
. Then%
\begin{equation*}
\clm\mathcal{F}\left( \left( \bar{A},\overline{b}\right) ,\overline{x}%
\right) =\left( \left\Vert \overline{x}\right\Vert +1\right) \clm\mathcal{F}%
_{\bar{A}}\left( \overline{b},\overline{x}\right)
\end{equation*}%
and
\begin{eqnarray*}
\clm\mathcal{F}_{\bar{A}}\left( \overline{b},\overline{x}\right)
&=&\max_{D\in \mathcal{D}_{\bar{A},\overline{b}}\left( x\right) }\left[ \dist%
\,_{\ast }\left( 0_{n},\conv\left\{ \overline{a}_{t},~t\in D\right\} \right) %
\right] ^{-1} \\
&=&\dist\,_{\ast }\left( 0_{n},\mathrm{end}\conv\left\{ \overline{a}%
_{t},~t\in T_{\bar{A},\overline{b}}\left( \overline{x}\right) \right\}
\right) ^{-1}.
\end{eqnarray*}
\end{theo}

The reader is addressed to \cite[Corollary 2.1, Remark 2.3 and Corollary 3.2]%
{LMY18} for an extension of the previous theorem to semi-infinite systems
(with infinitely many constraints) under a certain regularity condition.
Regarding semi-infinite systems, the monograph \cite{libro} contains a
comprehensive development of theory and applications of these models. See
also \cite{KNT10} for an approach of the calmness modulus for inequality
systems in terms of the subdifferential of an associated supremum function.

\section{Lipschitz upper semicontinuity of the feasible set mapping under
full perturbations}

This section tackles the main goal of the paper, which is to establish the
relationship between the Lipschitz upper semicontinuity modulus and the
calmness moduli of the feasible set mapping under full perturbations,
closing the gap shown in \eqref{eq Lipusc vs clm}. In a second step, a
computable point-based expression for the Lipschitz upper semicontinuity
modulus is provided through the calmness moduli at certain points of the
domain. The section is divided into two since the case of unbounded feasible
set has a pathological behavior in the context of full perturbations.

The following proposition reformulates the expression of the Lipschitz upper
semicontinuity modulus given in (\ref{eq Lipusc equiv}) in the context of
full perturbations of the system's parameters.

\begin{prop}
\label{prop Lipusc reformulation} Let $\left( \bar{A},\overline{b}\right)
\in \dom\mathcal{F}$. Then
\begin{eqnarray}
&&\Lipusc\mathcal{F}\left( \bar{A},\overline{b}\right)   \label{eq_5} \\
&=&\lim_{\delta \downarrow 0}\left( \sup \left\{ \frac{\left( \left\Vert
x\right\Vert +1\right) \dist\left( x,\mathcal{F}\left( \bar{A},\overline{b}%
\right) \right) }{\left\Vert \left( \bar{A}x-\overline{b}\right)
_{+}\right\Vert _{\infty }}\Bigg|\left\Vert \left( \bar{A}x-\overline{b}%
\right) _{+}\right\Vert _{\infty }\leq \delta \left( \left\Vert x\right\Vert
+1\right) \right\} \right) .  \notag
\end{eqnarray}
\end{prop}

\begin{dem}
First, observe that the supremum in the right-hand side of (\ref{eq_5}) is a
(real-extended) decreasing function of $\delta >0,$ so that the limit of the
supremum exists in $\left[ 0,+\infty \right] .$One can check, according to
the definitions and appealing to Lemma \ref{lem dist inverse} that, when the
right-hand side of (\ref{eq_5}) is infinite, then $\Lipusc\mathcal{F}\left(
\bar{A},\overline{b}\right) =+\infty .$ Assume now that this right-hand side
is finite and consider a paired sequence $\left( \delta _{k},x_{k}\right)
\in \mathbb{R}\times \mathbb{R}^{n}$ for which the limit of the supremum is
attained. More precisely, consider a sequence $\{\delta _{k}\}_{k\in \mathbb{%
N}}\downarrow 0$ and a corresponding sequence $\{x_{k}\}_{k\in \mathbb{N}}$
satisfying conditions
\begin{equation}
\left\Vert \left( \bar{A}x_{k}-\overline{b}\right) _{+}\right\Vert _{\infty
}\leq \delta _{k}\left( \left\Vert x_{k}\right\Vert +1\right) ,
\label{eq delta}
\end{equation}%
and
\begin{eqnarray}
&&\frac{\left( \left\Vert x_{k}\right\Vert +1\right) \dist\left( x_{k},%
\mathcal{F}\left( \bar{A},\overline{b}\right) \right) }{\left\Vert \left(
\bar{A}x_{k}-\overline{b}\right) _{+}\right\Vert _{\infty }}
\label{eq 15bis} \\
&\geq &\sup \left\{ \frac{\left( \left\Vert x\right\Vert +1\right) \dist%
\left( x,\mathcal{F}\left( \bar{A},\overline{b}\right) \right) }{\left\Vert
\left( \bar{A}x-\overline{b}\right) _{+}\right\Vert _{\infty }}\Bigg|%
\left\Vert \left( \bar{A}x-\overline{b}\right) _{+}\right\Vert _{\infty
}\leq \delta _{k}\left( \left\Vert x\right\Vert +1\right) \right\} -\frac{1}{%
k}.  \notag
\end{eqnarray}%
For each point $x_{k}$, $\mathcal{F}^{-1}(x_{k})$ is a closed subset (a
polyhedral convex cone, indeed) of $\mathbb{R}^{m\times n}\times \mathbb{R}%
^{m},$ so there exists a pair of parameters $\left( A_{k},b_{k}\right) \in
\mathcal{F}^{-1}(x_{k})$ such that
\begin{equation*}
\dist\left( \left( \bar{A},\overline{b}\right) ,\mathcal{F}%
^{-1}(x_{k})\right) =\left\Vert \left( \bar{A}-A_{k},\overline{b}%
-b_{k}\right) \right\Vert ,
\end{equation*}%
and, as a consequence of Lemma \ref{lem dist inverse}, we have
\begin{equation}
\left\Vert \left( \bar{A}-A_{k},\overline{b}-b_{k}\right) \right\Vert =\frac{%
\left\Vert \left( \bar{A}x_{k}-\overline{b}\right) _{+}\right\Vert _{\infty }%
}{\left\Vert x_{k}\right\Vert +1}\leq \delta _{k}.  \label{eq dist param}
\end{equation}%
Thus, we have:
\begin{align*}
\Lipusc\mathcal{F}\left( \bar{A},\overline{b}\right) & \geq \underset{%
k\rightarrow \infty }{\lim \sup }\frac{\dist\left( x_{k},\mathcal{F}\left(
\bar{A},\overline{b}\right) \right) }{\left\Vert \left( \bar{A}-A_{k},%
\overline{b}-b_{k}\right) \right\Vert } \\
& =\underset{k\rightarrow \infty }{\lim \sup }\frac{\left( \left\Vert
x_{k}\right\Vert +1\right) \dist\left( x_{k},\mathcal{F}\left( \bar{A},%
\overline{b}\right) \right) }{\left\Vert \left( \bar{A}x_{k}-\overline{b}%
\right) _{+}\right\Vert _{\infty }},
\end{align*}%
where the inequality comes from \eqref{eq Lipusc equiv} and the equality
from \eqref{eq dist param}. Hence due to the construction of the sequences $%
\{\delta _{k}\}_{k\in \mathbb{N}}$ and $\{x_{k}\}_{k\in \mathbb{N}}$ we have
\begin{eqnarray*}
&&\Lipusc\mathcal{F}\left( \bar{A},\overline{b}\right)  \\
&\geq &\lim_{\delta \downarrow 0}\left( \sup \left\{ \frac{\left( \left\Vert
x\right\Vert +1\right) \dist\left( x,\mathcal{F}\left( \bar{A},\overline{b}%
\right) \right) }{\left\Vert \left( \bar{A}x-\overline{b}\right)
_{+}\right\Vert _{\infty }}\Bigg|\left\Vert \left( \bar{A}x-\overline{b}%
\right) _{+}\right\Vert _{\infty }\leq \delta \left( \left\Vert x\right\Vert
+1\right) \right\} \right) .
\end{eqnarray*}%
On the other hand, in order to prove the converse inequality, consider
sequences $\{(A_{k},b_{k})\}_{k\in \mathbb{N}}\rightarrow \left( \bar{A},%
\overline{b}\right) $ and $\{x_{k}\}_{k\in \mathbb{N}}$ with $x_{k}\in
\mathcal{F}(A_{k},b_{k})$ for each $k$ such that
\begin{equation*}
\Lipusc\mathcal{F}\left( \bar{A},\overline{b}\right) =\lim_{k\rightarrow
\infty }\frac{\dist\left( x_{k},\mathcal{F}\left( \bar{A},\overline{b}%
\right) \right) }{\left\Vert \left( \bar{A}-A_{k},\overline{b}-b_{k}\right)
\right\Vert }.
\end{equation*}%
In such a case, we may define the corresponding sequence $\{\delta
_{k}\}_{k\in \mathbb{N}}\rightarrow 0$ by
\begin{equation*}
\delta _{k}:=\left\Vert \left( \bar{A}-A_{k},\overline{b}-b_{k}\right)
\right\Vert \geq \dist\left( \left( \bar{A},\overline{b}\right) ,\mathcal{F}%
^{-1}(x_{k})\right) =\frac{\left\Vert \left( \bar{A}x_{k}-\overline{b}%
\right) _{+}\right\Vert _{\infty }}{\left\Vert x_{k}\right\Vert +1},
\end{equation*}%
appealing again to Lemma \ref{lem dist inverse}. Hence, we can conclude that
\begin{align*}
& \Lipusc\mathcal{F}\left( \bar{A},\overline{b}\right)  \\
& \leq \underset{k\rightarrow \infty }{\lim \sup }\frac{\left( \left\Vert
x_{k}\right\Vert +1\right) \dist\left( x_{k},\mathcal{F}\left( \bar{A},%
\overline{b}\right) \right) }{\left\Vert \left( \bar{A}x_{k}-\overline{b}%
\right) _{+}\right\Vert _{\infty }} \\
& \leq \limsup_{k\rightarrow \infty }\left( \sup \left\{ \frac{\left(
\left\Vert x\right\Vert +1\right) \dist\left( x,\mathcal{F}\left( \bar{A},%
\overline{b}\right) \right) }{\left\Vert \left( \bar{A}x-\overline{b}\right)
_{+}\right\Vert _{\infty }}\Bigg|\left\Vert \left( \bar{A}x-\overline{b}%
\right) _{+}\right\Vert _{\infty }\leq \delta _{k}\left( \left\Vert
x\right\Vert +1\right) \right\} \right)  \\
& \leq \lim_{\delta \downarrow 0}\left( \sup \left\{ \frac{\left( \left\Vert
x\right\Vert +1\right) \dist\left( x,\mathcal{F}\left( \bar{A},\overline{b}%
\right) \right) }{\left\Vert \left( \bar{A}x-\overline{b}\right)
_{+}\right\Vert _{\infty }}\Bigg|\left\Vert \left( \bar{A}x-\overline{b}%
\right) _{+}\right\Vert _{\infty }\leq \delta \left( \left\Vert x\right\Vert
+1\right) \right\} \right) .
\end{align*}
\end{dem}

\subsection{Unbounded feasible set}

This subsection is devoted to the case where the nominal feasible set, $%
\mathcal{F}\left( \bar{A},\overline{b}\right) $, is unbounded. Indeed, the
Lipschitz upper semicontinuity modulus is infinite in this framework of full
perturbations, excluding a trivial case $\mathcal{F}\left( \bar{A},\overline{%
b}\right) =\mathbb{R}^{n}$ and the special case $n=1$, as opposed to the RHS
perturbations one. The following result formalizes this statement.

\begin{prop}
\label{Prop_F(A,b)_unbounded}Let $\left( \bar{A},\overline{b}\right) \in \dom%
\mathcal{F}$\thinspace\ with $\mathcal{F}\left( \bar{A},\overline{b}\right) $
unbounded. Then
\begin{equation*}
\Lipusc\mathcal{F}\left( \bar{A},\overline{b}\right) =\left\{
\begin{tabular}{ll}
$0$ &
\begin{tabular}{l}
$\text{if }\mathcal{F}\left( \bar{A},\overline{b}\right) =\mathbb{R}^{n},$%
\end{tabular}%
\medskip  \\
$\dfrac{\left\Vert \overline{x}\right\Vert +1}{\max \left\{ \left\Vert
\overline{a}_{t}\right\Vert _{\ast },~t\in T_{\bar{A},\overline{b}}\left(
\overline{x}\right) \right\} }$ & $%
\begin{array}{l}
\text{if }n=1\text{ and }\mathcal{F}\left( \bar{A},\overline{b}\right) \text{
is} \\
\left] -\infty ,\overline{x}\right] \text{ or }\left[ \overline{x},+\infty %
\right[ ,%
\end{array}%
$ \\
$+\infty $ &
\begin{tabular}{l}
$\text{otherwise}$.%
\end{tabular}%
\end{tabular}%
\right.
\end{equation*}
\end{prop}

\begin{dem}
The case $\mathcal{F}\left( \bar{A},\overline{b}\right) =\mathbb{R}^{n}$
follows trivially from the definition. If $n=1$ we will consider the case $%
\mathcal{F}\left( \bar{A},\overline{b}\right) =\left] -\infty ,\overline{x}%
\right] ,$ since the case $\left[ \overline{x},+\infty \right[ $ is
completely analogous. It is a well-known fact that in the case $n=1$ the
feasible set mapping $\mathcal{F}$ is upper semicontinuous in the sense of
Berge (see, e.g., \cite[Exercise 6.6]{libro}). This means that for all $%
\varepsilon >0$ there exists $\delta >0$ such that $\left\Vert \left(
A,b\right) -\left( \bar{A},\overline{b}\right) \right\Vert <\delta $ implies
$\mathcal{F}\left( A,b\right) \subset \left] -\infty ,\overline{x}%
+\varepsilon \right[ .$ Accordingly, looking at the definitions (\ref{eq
Lipusc property}) and (\ref{eq calmness property}), the Lipschitz upper
semicontinuity of $\mathcal{F}$ at $\left( \bar{A},\overline{b}\right) $
turns out to be equivalent to the calmness of $\mathcal{F}$ at $\left(
\left( \bar{A},\overline{b}\right) ,\overline{x}\right) $ and both moduli do
coincide. Since the case $\mathcal{F}\left( \bar{A},\overline{b}\right) =%
\left] -\infty ,\overline{x}\right] $ obviously implies $\overline{a}%
_{t}\geq 0$ for all $t\in T,$ one has
\begin{equation*}
\dist\,_{\ast }\left( 0_{n},\mathrm{end}\conv\left\{ \overline{a}_{t},~t\in
T_{\bar{A},\overline{b}}\left( \overline{x}\right) \right\} \right) =\max
\left\{ \left\Vert \overline{a}_{t}\right\Vert _{\ast },~t\in T_{\bar{A},%
\overline{b}}\left( \overline{x}\right) \right\} ,
\end{equation*}%
so that, recalling Theorem \ref{theo clm feasible}, we have
\begin{equation*}
\Lipusc\mathcal{F}\left( \bar{A},\overline{b}\right) =\clm\mathcal{F}\left(
\left( \bar{A},\overline{b}\right) ,\overline{x}\right) =\dfrac{\left\Vert
\overline{x}\right\Vert +1}{\max \left\{ \left\Vert \overline{a}%
_{t}\right\Vert _{\ast },~t\in T_{\bar{A},\overline{b}}\left( \overline{x}%
\right) \right\} }.
\end{equation*}

Assume now the situation $\mathcal{F}\left( \bar{A},\overline{b}\right)
\varsubsetneq \mathbb{R}^{n}$ with $n\geq 2.$ In such a case, every
unbounded closed convex set different from $\mathbb{R}^{n}$ is known to have
an unbounded boundary. So, consider an unbounded sequence $\{z_{k}\}_{k\in
\mathbb{N}}\subset \bd\mathcal{F}\left( \bar{A},\overline{b}\right) $. It is
not restrictive to assume $\left\Vert z_{k}\right\Vert \rightarrow +\infty .$
Then, for any $k\in \mathbb{N}$, $\clm\mathcal{F}_{\bar{A}}\left( \overline{b%
},z_{k}\right) >0$ by Theorem \ref{theo clm feasible}. More precisely, the
only possibility enabling $\clm\mathcal{F}_{\bar{A}}\left( \overline{b},%
\overline{x}\right) =0$ is $\mathrm{end}\conv\left\{ \overline{a}_{t},~t\in
T_{\bar{A},\overline{b}}\left( \overline{x}\right) \right\} =\emptyset ,$
and this never happens if $\overline{x}\in \bd\mathcal{F}\left( \bar{A},%
\overline{b}\right) $ because the separation theorem ensures the existence
of $t\in T_{\bar{A},\overline{b}}\left( \overline{x}\right) $ with $%
\overline{a}_{t}\neq 0_{n}.$ Since there are only finitely many subsets of $T
$, there exists some constant $\gamma >0$ such that $\clm\mathcal{F}_{\bar{A}%
}\left( \overline{b},z_{k}\right) \geq \gamma $ for all $k\in \mathbb{N}$.
By definition, we can find a sequence $\{x_{k}\}_{k\in \mathbb{N}}$
fulfilling $\Vert x_{k}-z_{k}\Vert \leq 1/k$ and
\begin{equation*}
\dist\left( x_{k},\mathcal{F}\left( \bar{A},\overline{b}\right) \right) \geq
\left( \clm\mathcal{F}_{\bar{A}}\left( \overline{b},z_{k}\right)
-(1/k)\right) \left\Vert \left( \bar{A}x_{k}-\overline{b}\right)
_{+}\right\Vert _{\infty }.
\end{equation*}%
Clearly, we have $\Vert x_{k}\Vert \rightarrow +\infty $ and $\left\Vert
\left( \bar{A}x_{k}-\overline{b}\right) _{+}\right\Vert _{\infty
}\rightarrow 0$, so that
\begin{equation*}
\left\Vert \left( \bar{A}x_{k}-\overline{b}\right) _{+}\right\Vert _{\infty
}/\left( \Vert x_{k}\Vert +1\right) \rightarrow 0.
\end{equation*}%
Taking into account Proposition \ref{prop Lipusc reformulation}, we conclude
that
\begin{align*}
\Lipusc\mathcal{F}\left( \bar{A},\overline{b}\right) & \geq
\limsup_{k\rightarrow \infty }\frac{\left( \left\Vert x_{k}\right\Vert
+1\right) \dist\left( x_{k},\mathcal{F}\left( \bar{A},\overline{b}\right)
\right) }{\left\Vert \left( \bar{A}x_{k}-\overline{b}\right) _{+}\right\Vert
_{\infty }} \\
& \geq \limsup_{k\rightarrow \infty }\left( \left\Vert x_{k}\right\Vert
+1\right) \left( \clm\mathcal{F}_{\bar{A}}\left( \overline{b},z_{k}\right)
-(1/k)\right)  \\
& \geq \limsup_{k\rightarrow \infty }\left( \left\Vert x_{k}\right\Vert
+1\right) \left( \gamma -(1/k)\right) =+\infty .
\end{align*}
\end{dem}

\subsection{Bounded feasible set}

This subsection deals with systems whose feasible set is bounded at the
nominal parameters. The main results are Theorem \ref{theo sup clm}, which
provides the equality between the Lipschitz upper semicontinuity modulus and
the supremum of calmness moduli at nominal feasible points, as well as
Theorem \ref{theo max clm}, which replaces the supremum with a maximum over
the set of extreme points of the feasible set. To achieve the second goal,
the technical Lemma \ref{lem tech clm extr} exploits the polyhedral
structure of the feasible set as an essential tool.

\begin{theo}
\label{theo sup clm} Let $\left( \bar{A},\overline{b}\right) \in \dom%
\mathcal{F}$ such that $\mathcal{F}\left( \bar{A},\overline{b}\right) $ is
bounded. Then
\begin{equation}
\Lipusc\mathcal{F}\left( \bar{A},\overline{b}\right) =\sup_{x\in \mathcal{F}%
\left( \bar{A},\overline{b}\right) }\clm\mathcal{F}\left( \left( \bar{A},%
\overline{b}\right) ,x\right) .  \label{Lipusc=sup full}
\end{equation}
\end{theo}

\begin{dem}
According to (\ref{eq Lipusc vs clm}), we only have to prove inequality $%
\leq .$ Appealing to Proposition \ref{prop Lipusc reformulation}, there
exists a paired sequence $\left( \delta _{k},x_{k}\right) _{k\in \mathbb{N}}$
satisfying $\delta _{k}\downarrow 0$ as well as (\ref{eq delta}) and (\ref%
{eq 15bis}) such that
\begin{equation*}
\Lipusc\mathcal{F}\left( \bar{A},\overline{b}\right) =\lim_{k\rightarrow
\infty }\frac{\left( \left\Vert x_{k}\right\Vert +1\right) \dist\left( x_{k},%
\mathcal{F}\left( \bar{A},\overline{b}\right) \right) }{\left\Vert \left(
\bar{A}x_{k}-\overline{b}\right) _{+}\right\Vert _{\infty }}.
\end{equation*}%
From (\ref{eq_Hof_lip}) and Remark \ref{Rem 1b} (see alternatively \cite[%
p.263]{Hof52}), there exists some constant $\kappa \geq 0$ fulfilling
\begin{equation}
\Vert x_{k}\Vert -\beta \leq \dist\left( x_{k},\mathcal{F}\left( \bar{A},%
\overline{b}\right) \right) \leq \kappa \left\Vert \left( \bar{A}x_{k}-%
\overline{b}\right) _{+}\right\Vert _{\infty }\leq \kappa \delta _{k}\left(
\left\Vert x_{k}\right\Vert +1\right) ,  \label{eq Hoffman bounds}
\end{equation}%
where $\beta :=\max \{\Vert x\Vert \mid x\in \mathcal{F}\left( \bar{A},%
\overline{b}\right) \}$. Thus, $\Vert x_{k}\Vert (1-\kappa \delta _{k})\leq
\beta +\kappa \delta _{k}$ and we can conclude that the sequence $%
\{x_{k}\}_{k\in \mathbb{N}}$ is bounded. Then, it is clear that $\kappa
\delta _{k}\left( \left\Vert x_{k}\right\Vert +1\right) \rightarrow 0$ and,
as a consequence of inequality \eqref{eq
Hoffman bounds}, we have $\dist\left( x_{k},\mathcal{F}\left( \bar{A},%
\overline{b}\right) \right) \rightarrow 0$ too. Now consider for every $k\in
\mathbb{N}$ a projection of $x_{k}$ onto $\mathcal{F}\left( \bar{A},%
\overline{b}\right) $, denoted by $z_{k}\in \mathcal{F}\left( \bar{A},%
\overline{b}\right) $. In other words, $\Vert x_{k}-z_{k}\Vert =\dist\left(
x_{k},\mathcal{F}\left( \bar{A},\overline{b}\right) \right) \rightarrow 0$.
Then, we also have
\begin{equation}
\Lipusc\mathcal{F}\left( \bar{A},\overline{b}\right) =\lim_{k\rightarrow
\infty }\frac{\left( \left\Vert z_{k}\right\Vert +1\right) \dist\left( x_{k},%
\mathcal{F}\left( \bar{A},\overline{b}\right) \right) }{\left\Vert \left(
\bar{A}x_{k}-\overline{b}\right) _{+}\right\Vert _{\infty }}.
\label{eq Lipusc zk}
\end{equation}%
We may consider any convex combination between $x_{k}$ and $z_{k}$, say $%
x_{k}^{\lambda }:=\lambda x_{k}+(1-\lambda )z_{k}$ for $\lambda \in \lbrack
0,1]$. Indeed, $z_{k}$ is still a projection of $x_{k}^{\lambda }$ onto $%
\mathcal{F}\left( \bar{A},\overline{b}\right) $ and
\begin{equation*}
\dist\left( x_{k}^{\lambda },\mathcal{F}\left( \bar{A},\overline{b}\right)
\right) =\lambda \dist\left( x_{k},\mathcal{F}\left( \bar{A},\overline{b}%
\right) \right) .
\end{equation*}%
Moreover, we have $\Vert (\bar{A}x_{k}^{\lambda }-\overline{b})_{+}\Vert
_{\infty }\leq \lambda \Vert (\bar{A}x_{k}-\overline{b})_{+}\Vert _{\infty }$%
, implying
\begin{equation*}
\frac{\dist\left( x_{k},\mathcal{F}\left( \bar{A},\overline{b}\right)
\right) }{\left\Vert \left( \bar{A}x_{k}-\overline{b}\right) _{+}\right\Vert
_{\infty }}\leq \frac{\dist\left( x_{k}^{\lambda },\mathcal{F}\left( \bar{A},%
\overline{b}\right) \right) }{\left\Vert \left( \bar{A}x_{k}^{\lambda }-%
\overline{b}\right) _{+}\right\Vert _{\infty }}\quad \forall \lambda \in
\lbrack 0,1].
\end{equation*}%
Passing to the limit on $\lambda \rightarrow 0$, we obtain
\begin{equation*}
\frac{\dist\left( x_{k},\mathcal{F}\left( \bar{A},\overline{b}\right)
\right) }{\left\Vert \left( \bar{A}x_{k}-\overline{b}\right) _{+}\right\Vert
_{\infty }}\leq \limsup_{\lambda \rightarrow 0}\frac{\dist\left(
x_{k}^{\lambda },\mathcal{F}\left( \bar{A},\overline{b}\right) \right) }{%
\left\Vert \left( \bar{A}x_{k}^{\lambda }-\overline{b}\right)
_{+}\right\Vert _{\infty }}\leq \clm\mathcal{F}\left( \left( \bar{A},%
\overline{b}\right) ,z_{k}\right) .
\end{equation*}%
Finally, combining the last inequality with expression \eqref{eq Lipusc zk},
we conclude that
\begin{equation*}
\Lipusc\mathcal{F}\left( \bar{A},\overline{b}\right) \leq \sup_{x\in
\mathcal{F}\left( \bar{A},\overline{b}\right) }\clm\mathcal{F}\left( \left(
\bar{A},\overline{b}\right) ,x\right) ,
\end{equation*}%
as desired.
\end{dem}

\begin{rem}
\emph{Indeed, formula (\ref{Lipusc=sup full}) holds true for every }$\left(
\bar{A},\overline{b}\right) \in \dom\mathcal{F}.$\emph{\ More specifically,
in the case} $\mathcal{F}\left( \bar{A},\overline{b}\right) =\mathbb{R}^{n}$%
\emph{\ both sides in (\ref{Lipusc=sup full}) are trivially zero from the
definitions; in the case }$n=1$\emph{\ equality follows from the previous
theorem and Proposition \ref{Prop_F(A,b)_unbounded}; and the case }$n\geq 2$%
\emph{\ with }$\mathcal{F}\left( \bar{A},\overline{b}\right) $\emph{\
unbounded and different from }$\mathbb{R}^{n}$\emph{\ we may proceed as in
the proof of Proposition \ref{Prop_F(A,b)_unbounded}, finding a sequence }$%
\{z_{k}\}_{k\in \mathbb{N}}\subset \bd\mathcal{F}\left( \bar{A},\overline{b}%
\right) $\emph{\ with }$\left\Vert z_{k}\right\Vert \rightarrow +\infty $%
\emph{\ and considering }$\gamma >0$\emph{\ with }$\clm\mathcal{F}_{\bar{A}%
}\left( \overline{b},z_{k}\right) \geq \gamma $\emph{\ for all }$k.$ \emph{%
Then }$\clm\mathcal{F}\left( \left( \bar{A},\overline{b}\right)
,z_{k}\right) \geq \left( \left\Vert z_{k}\right\Vert +1\right) \gamma
\rightarrow +\infty .$\emph{\ Accordingly, both sides of (\ref{Lipusc=sup
full}) are infinite.}
\end{rem}

This result completes our first goal, which was to find the counterpart of
Theorem \ref{theo Lipusc vs clm convex} in our framework. The following
lemma relies on the polyhedral structure of the feasible set mapping $%
\mathcal{F}_{\bar{A}}$ to provide some hierarchy among the calmness moduli
at different points of the image set. It follows the same spirit of the
proof of \cite[Proposition 5 and Theorem 6]{CCP21}, but coupled here in a
single result for completeness.

\begin{lem}
\label{lem tech clm extr} Let $F$ be a face of $\mathcal{F}\left( \bar{A},%
\overline{b}\right) $ for $\left( \bar{A},\overline{b}\right) \in \dom%
\mathcal{F}$. Then, $\clm\mathcal{F}_{\bar{A}}\left( \overline{b},\cdot
\right) $ is constant on $\ri F$, i.e., there exists $\gamma _{F}\geq 0$
such that $\clm\mathcal{F}_{\bar{A}}\left( \overline{b},x\right) =\gamma
_{F} $ for all $x\in \ri F$. Further,
\begin{equation*}
\clm\mathcal{F}_{\bar{A}}\left( \overline{b},x\right) \geq \gamma _{F}\quad
\text{for all }x\in F\setminus \ri F.
\end{equation*}
\end{lem}

\begin{dem}
Take any $x\in \ri F$ and set $I:=T_{\bar{A},\overline{b}}(x)$. Consider the
associated set
\begin{equation*}
F_{I}:=\left\{ \tilde{x}\in F\Bigg|%
\begin{tabular}{rl}
$\overline{a}_{t}^{\prime }\tilde{x}=\overline{b}_{t}$, & $t\in I$, \\
$\overline{a}_{t}^{\prime }\tilde{x}\leq \overline{b}_{t}$, & $t\in
T\setminus I$.%
\end{tabular}%
\ \right\} .
\end{equation*}%
By construction, $F_{I}$ is a face of $F$ and $x\in \ri F_{I}$ (see, e.g.,
\cite[Theorem A.7]{libro}). Hence, appealing to \cite[Corollary 18.1.2]{Ro70}%
, this means that $F=F_{I}$ and so
\begin{equation*}
\ri F=\left\{ x\in F\Bigg|%
\begin{tabular}{rl}
$\overline{a}_{t}^{\prime }x=\overline{b}_{t}$, & $t\in I$, \\
$\overline{a}_{t}^{\prime }x<\overline{b}_{t}$, & $t\in T\setminus I$%
\end{tabular}%
\ \right\} .
\end{equation*}%
Therefore, $T_{\bar{A},\overline{b}}(x)=I$ for all $x\in \ri F$ and, as a
consequence of Theorem \ref{theo clm feasible}, we conclude that $\clm%
\mathcal{F}_{\bar{A}}\left( \overline{b},\cdot \right) $ is constant on $\ri %
F$. For the second part of the statement, consider any $\hat{x}\in
F\setminus \ri F$ and any index set $D\in \mathcal{D}_{\bar{A},\overline{b}%
}(x)$ with an associated solution, $d\in \mathbb{R}^{n}$, of the system %
\eqref{eq system D}. Since $D\subset I\subset T_{\bar{A},\overline{b}}(\hat{x%
})$, we obtain that
\begin{equation*}
\left\{
\begin{tabular}{rl}
$\overline{a}_{t}^{\prime }\left( x-\hat{x}\right) =0$, & $t\in D$, \\
$\overline{a}_{t}^{\prime }\left( x-\hat{x}\right) <0$, & $t\in T_{\bar{A},%
\overline{b}}(\hat{x})\setminus D$%
\end{tabular}%
\ \right\} ,
\end{equation*}%
and consequently
\begin{equation*}
\left\{
\begin{tabular}{rl}
$\overline{a}_{t}^{\prime }\left( d+\lambda \left( x-\hat{x}\right) \right) =%
\overline{a}_{t}^{\prime }d=1$, & $t\in D$, \\
$\overline{a}_{t}^{\prime }\left( d+\lambda \left( x-\hat{x}\right) \right) =%
\overline{a}_{t}^{\prime }d<1$, & $t\in I\setminus D$, \\
$\overline{a}_{t}^{\prime }\left( d+\lambda \left( x-\hat{x}\right) \right) <%
\overline{a}_{t}^{\prime }d<1$, & $t\in T_{\bar{A},\overline{b}}(\hat{x}%
)\setminus I$%
\end{tabular}%
\ \right\} ,
\end{equation*}%
for any $\lambda >0$. Thus, $D\in \mathcal{D}_{\bar{A},\overline{b}}(\hat{x})
$ and $\clm\mathcal{F}_{\bar{A}}\left( \overline{b},\hat{x}\right) \geq \clm%
\mathcal{F}_{\bar{A}}\left( \overline{b},x\right) $ follows from Theorem \ref%
{theo clm feasible}.
\end{dem}

\begin{theo}
\label{theo max clm} Let $\left( \bar{A},\overline{b}\right) \in \dom%
\mathcal{F}$ with $\mathcal{F}\left( \bar{A},\overline{b}\right) $ bounded.
Then
\begin{equation*}
\Lipusc\mathcal{F}\left( \bar{A},\overline{b}\right) =\max_{x\in \mathcal{E}%
\left( \bar{A},\overline{b}\right) }\clm\mathcal{F}\left( \left( \bar{A},%
\overline{b}\right) ,x\right) .
\end{equation*}
\end{theo}

\begin{dem}
By Theorem \ref{theo sup clm}, we have
\begin{equation*}
\Lipusc\mathcal{F}\left( \bar{A},\overline{b}\right) =\sup_{x\in \mathcal{F}%
\left( \bar{A},\overline{b}\right) }\clm\mathcal{F}\left( \left( \bar{A},%
\overline{b}\right) ,x\right) ,
\end{equation*}%
so it is enough to show that the supremum is attained at some extreme point
of the feasible set. Indeed, the supremum is known to be a maximum because of
Theorem \ref{theo clm feasible} and the finiteness of $T$ (or,
alternatively, from Lemma \ref{lem tech clm extr} and the finiteness of the
number of faces of $\mathcal{F}\left( \bar{A},\overline{b}\right) ;$ see
\cite[Theorem 19.1]{Ro70}). Write
\begin{equation*}
\Lipusc\mathcal{F}\left( \bar{A},\overline{b}\right) =\clm\mathcal{F}\left(
\left( \bar{A},\overline{b}\right) ,\overline{x}\right) .
\end{equation*}%
Then, setting $I:=T_{\bar{A},\overline{b}}(\overline{x}),$ we can construct
the face $F_{I}:=\{x\in \mathbb{R}^{n}\mid \overline{a}_{t}^{\prime }x=%
\overline{b}_{t},~t\in I;~\overline{a}_{t}^{\prime }x\leq \overline{b}%
_{t},~t\in T\setminus I\}$ of $\mathcal{F}\left( \bar{A},\overline{b}\right)
$ in which $\overline{x}$ belongs to its relative interior. Now, without
loss of generality (cf. \cite[Corollary 32.3.2]{Ro70}), take any extreme
point $\overline{z}$ on $F_{I}$ satisfying $\Vert \overline{z}\Vert =\max
\{\Vert z\Vert \mid z\in F_{I}\}$. According to Lemma \ref{lem tech clm extr}%
, we have
\begin{equation*}
\clm\mathcal{F}_{\bar{A}}\left( \overline{b},\overline{x}\right) \leq \clm%
\mathcal{F}_{\bar{A}}\left( \overline{b},\overline{z}\right) ,
\end{equation*}%
which implies%
\begin{equation*}
\left( \left\Vert \overline{x}\right\Vert +1\right) \clm\mathcal{F}_{\bar{A}%
}\left( \overline{b},\overline{x}\right) \leq \left( \left\Vert \overline{z}%
\right\Vert +1\right) \clm\mathcal{F}_{\bar{A}}\left( \overline{b},\overline{%
z}\right) .
\end{equation*}%
Observe that $\overline{z}$ is also an extreme point of $\mathcal{F}\left(
\bar{A},\overline{b}\right) $, cf. \cite[p.163]{Ro70}, so that by Theorem \ref{theo clm feasible} we obtain that
\begin{align*}
\Lipusc\mathcal{F}\left( \bar{A},\overline{b}\right) &=\clm\mathcal{F}\left(
\left( \bar{A},\overline{b}\right) ,\overline{x}\right) =\left( \Vert \overline{%
x}\Vert +1\right) \clm\mathcal{F}_{\bar{A}}\left( \overline{b},\overline{x}%
\right)\\
&\leq \left( \Vert \overline{%
z}\Vert +1\right) \clm\mathcal{F}_{\bar{A}}\left( \overline{b},\overline{z}%
\right)=\clm%
\mathcal{F}\left( \left( \bar{A},\overline{b}\right) ,\overline{z}\right)\\
&\leq \max_{x\in \mathcal{E}\left( \bar{A},\overline{b}\right) }\clm%
\mathcal{F}\left( \left( \bar{A},\overline{b}\right) ,x\right)
\leq \Lipusc\mathcal{F}\left( \bar{A},\overline{b}\right)
\end{align*}
and the assertion follows.
\end{dem}

We finish this section by providing a computable expression for the
Lipschitz upper semicontinuity modulus in our context of full perturbations.
The following corollary is a direct consequence of Theorems \ref{theo clm
feasible} and \ref{theo max clm}.

\begin{cor}
\label{Cor_Lipusc_F}Let $\left( \bar{A},\overline{b}\right) \in \dom\mathcal{%
F}$ with $\mathcal{F}\left( \bar{A},\overline{b}\right) $ bounded.\textbf{\ }%
Then
\begin{equation*}
\Lipusc\mathcal{F}\left( \bar{A},\overline{b}\right) =\max_{x\in \mathcal{E}%
\left( \bar{A},\overline{b}\right) }\left( \left( \left\Vert x\right\Vert
+1\right) \max_{D\in \mathcal{D}_{\bar{A},\overline{b}}\left( x\right) }%
\left[ \dist\,_{\ast }\left( 0_{n},\conv\left\{ \overline{a}_{t},~t\in
D\right\} \right) \right] ^{-1}\right) .
\end{equation*}
\end{cor}

\section{Conclusions, illustrations and further research}

For comparative purposes, let us gather here the expressions for the
Lipschitz upper semicontinuity moduli of $\mathcal{F}$ and $\mathcal{F}_{%
\bar{A}}.$ In the case $n=1$ we have that $\Lipusc\mathcal{F}\left( \bar{A},%
\overline{b}\right) $ is finite at any $\left( \bar{A},\overline{b}\right)
\in \dom\mathcal{F}.$ This case constitutes an exception,
since $\Lipusc\mathcal{F}\left( \bar{A},\overline{b}\right) =+\infty $ when $%
n\geq 2,$ $\mathcal{F}\left( \bar{A},\overline{b}\right) \varsubsetneq
\mathbb{R}^{n}$ and it is unbounded. Let us confine here to the case when $%
\mathcal{F}\left( \bar{A},\overline{b}\right) $ is a nonempty bounded set,
in which case $\mathcal{E}\left( \bar{A},\overline{b}\right) =\extr\mathcal{F%
}\left( \bar{A},\overline{b}\right) .$ Therefore,
\begin{equation}
\Lipusc\mathcal{F}_{\bar{A}}\left( \overline{b}\right) =\max_{x\in \extr%
\mathcal{F}\left( \bar{A},\overline{b}\right) }\clm\mathcal{F}_{\bar{A}%
}\left( \overline{b},x\right) ,  \label{eq_11}
\end{equation}%
while
\begin{eqnarray}
\Lipusc\mathcal{F}\left( \bar{A},\overline{b}\right)  &=&\max_{x\in \extr%
\mathcal{F}\left( \bar{A},\overline{b}\right) }\clm\mathcal{F}\left( \left(
\bar{A},\overline{b}\right) ,x\right)   \notag \\
&=&\max_{x\in \extr\mathcal{F}\left( \bar{A},\overline{b}\right) }\left(
\left\Vert x\right\Vert +1\right) \clm\mathcal{F}_{\bar{A}}\left( \overline{b%
},x\right) .  \label{eq_22}
\end{eqnarray}%
Moreover, it was previously known that, for any $\left( b,x\right) \in
\mathrm{gph}\mathcal{F}_{\bar{A}},$
\begin{equation}
\clm\mathcal{F}_{\bar{A}}\left( b,x\right) =\dist\,_{\ast }\left( 0_{n},%
\mathrm{end}\conv\left\{ \overline{a}_{t},~t\in T_{\bar{A},b}\left( x\right)
\right\} \right) ^{-1}.  \label{eq_00}
\end{equation}%
It is clear that the smaller the modulus, the greater the stability. The
following example shows, on the one hand, the dependence of these moduli on
the representation of the feasible set $\mathcal{F}\left( \bar{A},\overline{b%
}\right) $ (in particular, on the size of the coefficients of $\bar{A}$ and $%
\overline{b}).$ On the other hand, it is intended to illustrate the fact
that the Lipschitz upper semicontinuity modulus may not be a robust measure
in the sense that small perturbations of $\left( \bar{A},\overline{b}\right)
$ may cause great differences in $\Lipusc\mathcal{F}\left( \bar{A},\overline{%
b}\right) $, formally that function $\left( A,b\right) \mapsto \Lipusc%
\mathcal{F}\left( A,b\right) $ is not continuous at all $\left( \bar{A},%
\overline{b}\right) \in \mathrm{dom}\mathcal{F}$. Hence, the study of those $%
\left( \bar{A},\overline{b}\right) \in \mathrm{dom}\mathcal{F}$ where such a
function is continuous could constitute stuff for future research; see \cite%
{CCLP23} for counterpart results related to the continuity of the calmness
modulus.

\begin{exa}
\emph{\ Assume }$\mathbb{R}^{2}$ \emph{endowed with the Euclidean norm}
\emph{and} \emph{for each }$r\in \mathbb{N},$ \emph{let us consider the
system,}
\begin{equation*}
\sigma \left( A_{r},\overline{b}\right) :=\left\{ x_{1}\leq 1,\text{ }\frac{1%
}{r}x_{1}+\frac{1}{r^{2}}x_{2}\leq \frac{1}{r}-\frac{1}{r^{2}},\text{ }%
x_{2}\leq 0,\text{ }-x_{1}-x_{2}\leq 0\right\} .
\end{equation*}%
\emph{Observe that} $\left\{ \left( A_{r},\overline{b}\right) \right\}
_{r\in \mathbb{N}}$ \emph{converge to} $\left( \bar{A},\overline{b}\right) $
\emph{whose associated system is}
\begin{equation*}
\sigma \left( \bar{A},\overline{b}\right) :=\left\{ x_{1}\leq 1,\text{ }%
0x_{1}\leq 0,\text{ }x_{2}\leq 0,\text{ }-x_{1}-x_{2}\leq 0\right\} .
\end{equation*}

\emph{It is obvious that }$\mathcal{F}\left( A_{r},\overline{b}\right) =%
\mathrm{conv\{}\left( 0,0\right) ^{\prime },\left( 1-\frac{1}{r},0\right)
^{\prime },\left( 1,-1\right) ^{\prime }\},$ \emph{for all} $r\in \mathbb{N},
$ \emph{and }$\mathcal{F}\left( \bar{A},\overline{b}\right) =\mathrm{conv\{}%
\left( 0,0\right) ^{\prime },\left( 1,0\right) ^{\prime },\left( 1,-1\right)
^{\prime }\}.$ \emph{Applying (\ref{eq_11}), (\ref{eq_22}), and (\ref{eq_00}%
), and assuming }$r\geq 3,$\emph{\ one easily checks the following:}%
\begin{equation*}
\begin{tabular}{|l|l|l|l|}
\hline
$x\in \extr\mathcal{F}\left( A_{r},\overline{b}\right) $ & $T_{A_{r},%
\overline{b}}\left( x\right) $ & $\clm\mathcal{F}_{A_{r}}\left( \overline{b}%
,x\right) $ & $\clm\mathcal{F}\left( (A_{r},\overline{b}),x\right) $ \\
\hline
$\left( 0,0\right) ^{\prime }$ & $\{3,4\}$ & $\sqrt{5}$ & $\sqrt{5}$ \\
\hline
$\left( 1-\frac{1}{r},0\right) ^{\prime }$ & $\{2,3\}$ & $\sqrt{r^{2}-\frac{%
r^{2}-1}{r^{2}}}$ & $\left( 2-\frac{1}{r}\right) \sqrt{r^{2}-\frac{r^{2}-1}{%
r^{2}}}$ \\ \hline
$\left( 1,-1\right) ^{\prime }$ & $\{1,2,4\}$ & $\sqrt{r^{2}-\frac{r^{2}}{%
r^{2}+1}}$ & $\left( \sqrt{2}+1\right) \sqrt{r^{2}-\frac{r^{2}}{r^{2}+1}}$
\\ \hline
\end{tabular}%
\end{equation*}

\emph{Consequently,}%
\begin{eqnarray*}
\Lipusc\mathcal{F}_{A_{r}}\left( \overline{b}\right)  &=&\clm\mathcal{F}%
_{A_{r}}\left( \overline{b},\left( 1-\frac{1}{r},0\right) ^{\prime }\right) =%
\sqrt{r^{2}-\frac{r^{2}-1}{r^{2}}},\text{ } \\
\Lipusc\mathcal{F}\left( A_{r},\overline{b}\right)  &=&\clm\mathcal{F}\left(
(A_{r},\overline{b}),\left( 1,-1\right) ^{\prime }\right) =\left( \sqrt{2}%
+1\right) \sqrt{r^{2}-\frac{r^{2}}{r^{2}+1}}.
\end{eqnarray*}%
\emph{Thus }%
\begin{equation*}
\lim_{r\rightarrow \infty }\Lipusc\mathcal{F}_{A_{r}}\left( \overline{b}%
\right) =\lim_{r\rightarrow \infty }\Lipusc\mathcal{F}\left( A_{r},\overline{%
b}\right) =+\infty .
\end{equation*}%
\emph{However,} \emph{the situation regarding the limit system }$\sigma
\left( \bar{A},\overline{b}\right) $ \emph{is notably different:}
\begin{equation*}
\begin{tabular}{|l|l|l|l|}
\hline
$x\in \extr\mathcal{F}\left( \bar{A},\overline{b}\right) $ & $T_{\bar{A},%
\overline{b}}\left( x\right) $ & $\clm\mathcal{F}_{\bar{A}}\left( \overline{b%
},x\right) $ & $\clm\mathcal{F}\left( (\bar{A},\overline{b}),x\right) $ \\
\hline
$\left( 0,0\right) ^{\prime }$ & $\{3,4\}$ & $\sqrt{5}$ & $\sqrt{5}$ \\
\hline
$\left( 1,0\right) ^{\prime }$ & $\{1,2,3\}$ & $\sqrt{2}$ & $2\sqrt{2}$ \\
\hline
$\left( 1,-1\right) ^{\prime }$ & $\{1,2,4\}$ & $\sqrt{5}$ & $\left( \sqrt{2}%
+1\right) \sqrt{5}$ \\ \hline
\end{tabular}%
\end{equation*}%
\emph{So,}%
\begin{eqnarray*}
\Lipusc\mathcal{F}\left( \overline{b}\right)  &=&\clm\mathcal{F}_{\bar{A}%
}\left( \overline{b},\left( 0,0\right) ^{\prime }\right) =\clm\mathcal{F}_{%
\bar{A}}\left( \overline{b},\left( 1,-1\right) ^{\prime }\right) =\sqrt{5},%
\text{ } \\
\Lipusc\mathcal{F}\left( \bar{A},\overline{b}\right)  &=&\clm\mathcal{F}%
\left( (\bar{A},\overline{b}),\left( 1,-1\right) ^{\prime }\right) =\left(
\sqrt{2}+1\right) \sqrt{5}.
\end{eqnarray*}

\emph{The following picture is intended to illustrate the behavior of the
end sets} $\mathrm{end}\conv\left\{ a_{r,t},~t\in T_{A_{r},\overline{b}%
}\left( \overline{x}\right) \right\} $, \emph{where }$a_{r,t}$\emph{\ is the
}$t$\emph{-th row of }$A_{r},$ $r\in \mathbb{N},$ \emph{and} $\mathrm{end}%
\conv\left\{ \overline{a}_{t},~t\in T_{\bar{A},\overline{b}}\left( \overline{%
x}\right) \right\} $ \emph{at point} $\overline{x}=\left( 1,-1\right)
^{\prime }$ \emph{which is the one where we get the maximum of calmness
moduli }$\clm\mathcal{F}\left( (A_{r},\overline{b}),x\right) $ \emph{and} $%
\clm\mathcal{F}\left( (\bar{A},\overline{b}),x\right) $\emph{, among the
extreme points of the corresponding feasible sets. Looking at (\ref{eq_00})
the inverse of the distance from }$0_{2}$ \emph{to these end sets provides
the exact value of the calmness moduli and, consequently, in the computation
of the Lipschitz upper semicontinuity moduli.}

\begin{figure}[htbp]
    \centering
    \begin{minipage}{0.5\textwidth}
        \centering
        \includegraphics[width=\linewidth]{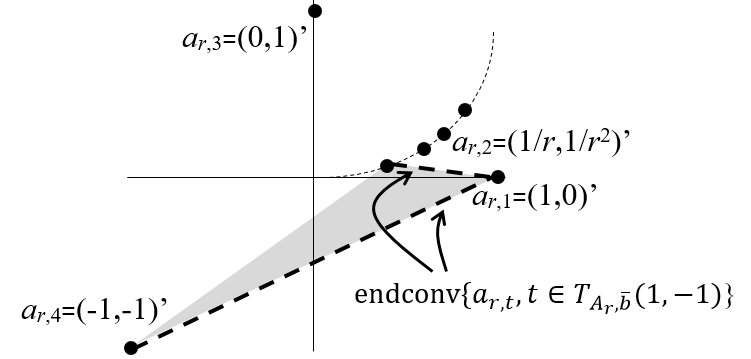}
    \end{minipage}\hfill
    \begin{minipage}{0.45\textwidth}
        \centering
        \includegraphics[width=\linewidth]{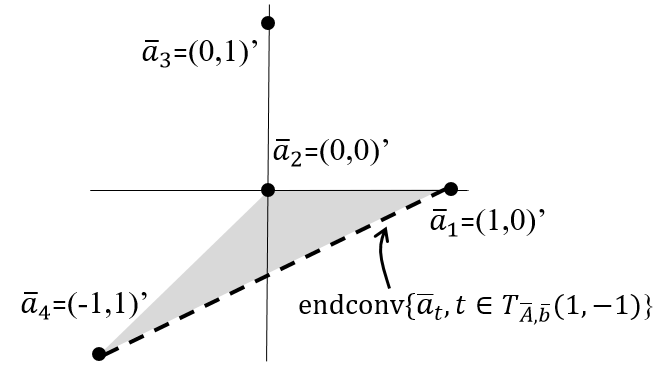}
    \end{minipage}
\end{figure}
\end{exa}

The following example shows that if we remove the first inequality in the
previous example, the situation is notably different

\begin{exa}
\emph{Assume again} $\mathbb{R}^{2}$ \emph{endowed with the Euclidean norm}
\emph{and,} \emph{for each }$r\in \mathbb{N},$ \emph{let us consider the
system,}
\begin{equation*}
\sigma \left( A_{r},\overline{b}\right) :=\left\{ \frac{1}{r}x_{1}+\frac{1}{%
r^{2}}x_{2}\leq \frac{1}{r}-\frac{1}{r^{2}},\text{ }x_{2}\leq 0,\text{ }%
-x_{1}-x_{2}\leq 0\right\} ,
\end{equation*}%
\begin{equation*}
\sigma \left( \bar{A},\overline{b}\right) :=\left\{ 0x_{1}\leq 0,\text{ }%
x_{2}\leq 0,\text{ }-x_{1}-x_{2}\leq 0\right\} .
\end{equation*}%
\emph{Here for }$r\geq 3$
\begin{gather*}
\Lipusc\mathcal{F}\left( A_{r},\overline{b}\right) =\clm\mathcal{F}\left(
(A_{r},\overline{b}),\left( 1,-1\right) ^{\prime }\right) =\left( \sqrt{2}%
+1\right) \sqrt{r^{2}-\frac{r^{2}}{r^{2}+1}}, \\
\clm\mathcal{F}\left( (\bar{A},\overline{b}),\left( 1,-1\right) ^{\prime
}\right) =\left( \sqrt{2}+1\right) \sqrt{2},
\end{gather*}%
\emph{however the unboundedness of }$\mathcal{F}(\bar{A},\overline{b})$\emph{%
\ yields }%
\begin{equation*}
\Lipusc\mathcal{F}\left( \bar{A},\overline{b}\right) =+\infty .
\end{equation*}
\end{exa}

\end{document}